\documentclass[11pt]{article}

\usepackage{amsfonts}
\usepackage{amsthm}
\usepackage{amsmath}
\usepackage{amssymb}
\usepackage{bm}
\usepackage{wasysym}

\usepackage{epsfig}
\usepackage[percent]{overpic}
\usepackage{graphics}
\usepackage{psfrag}
\usepackage[caption=false]{subfig}

\usepackage{titlesec}
\titleformat{\paragraph}
{\normalfont\normalsize\bfseries}{\theparagraph}{1em}{}
\titlespacing*{\paragraph}
{0pt}{3.25ex plus 1ex minus .2ex}{1.5ex plus .2ex}

\usepackage{verbatim}
\usepackage{color}
\usepackage{tikz}
\usetikzlibrary{arrows}
\usepackage{xcolor}
\definecolor{blun}{cmyk}{0.8, 0.5, 0, 0.7}

\usepackage{placeins}

\usepackage[all]{xy}
\usepackage{makeidx}
\usepackage{enumerate}
\usepackage[letterpaper,margin=0.9in]{geometry}
\usepackage[bookmarks,colorlinks,pdfauthor={DausKuehnSoresina},linkcolor=blun,citecolor=blun,urlcolor=blun]{hyperref}

\let\oldbibliography\thebibliography
\renewcommand{\thebibliography}[1]{%
  \oldbibliography{#1}%
  \setlength{\itemsep}{-1.2mm}%
}

\theoremstyle{plain}

\theoremstyle{definition}

\newtheoremstyle{myremark}
  {3pt}
  {3pt}
  {\small \rmfamily}
  {5pt}
  {\rmfamily}
  {:}
  {.5em}
  {}

\theoremstyle{myremark}

\newcommand{\be}{\begin{equation}}
\newcommand{\ee}{\end{equation}}
\newcommand{\benn}{\begin{equation*}}
\newcommand{\eenn}{\end{equation*}}
\newcommand{\bea}{\begin{eqnarray}}
\newcommand{\eea}{\end{eqnarray}}
\newcommand{\beann}{\begin{eqnarray*}}
\newcommand{\eeann}{\end{eqnarray*}}

\newcommand{\myendex}{$\blacklozenge$\end{ex}}
\newcommand{\myendexerc}{$\lozenge$\end{exerc}}
\newcommand{\myendpexerc}{$\lozenge$\end{pexerc}}

\begin{document}

\author{
Marzia Bisi\thanks{Department of Mathematical, Physical and Computer Sciences, University of Parma, Italy,  emails: \href{mailto: marzia.bisi@unipr.it}{marzia.bisi@unipr.it},  \href{mailto: maria.groppi@unipr.it}{maria.groppi@unipr.it}, \href{giorgio.martalo@unipr.it}{giorgio.martalo@unipr.it}}, 
Maria Groppi\footnotemark[1],
Giorgio Martal\`o\footnotemark[1],
Cinzia Soresina\thanks{Institute for Mathematics and Scientific Computing, University of Graz, 8010 Graz, Austria, email: \href{mailto: cinzia.soresina@uni-graz.at}{cinzia.soresina@uni-graz.at}}
}

\title{A chemotaxis reaction--diffusion model\\ for Multiple Sclerosis with Allee effect}

\maketitle

\begin{abstract} \noindent
In this paper, we study a modification of the mathematical model describing inflammation and demyelination patterns in the brain caused by Multiple Sclerosis proposed in [Lombardo et al. (2017), Journal of Mathematical Biology, 75, 373--417]. In particular, we hypothesize a minimal amount of macrophages to be able to start and sustain the inflammatory response. Thus, the model function for macrophage activation includes an Allee effect. We investigate the emergence of Turing patterns by combining linearised and weakly nonlinear analysis, bifurcation diagrams and numerical simulations, focusing on the comparison with the previous model.

{\begin{center} {\it  Dedicated to the memory of prof. Salvatore Rionero}\end{center}}



\end{abstract}

\section{Introduction}

Multiple Sclerosis is a demyelinating disorder affecting the central nervous system and causing severe and progressive physical and neurological impairment. More specifically, it is characterized by inflammation and demyelination, resulting in the formation of focal areas of myelin loss in the white matter of the brain, called plaques or lesions \cite{Lassmann}.
Different histological patterns of plaques have been identified \cite{Lucchinetti-etal}, but it is commonly believed that the onset of any typical new lesion is characterized by the same pathological changes \cite{Barnett-Prineas}.

Some mathematical models able to reproduce many of the typical pathological features of Multiple Sclerosis have been proposed and investigated; they are usually constituted by proper systems of PDEs \cite{Khonsari-Calvez-2007, Calvez-Khonsari-2008, barresi2016wavefront, lombardo2017demyelination}, but also systems of ODEs \cite{Frascoli-etal} or stochastic models \cite{Bordi-etal} are able to describe the relapse--remittance dynamics in patients affected by Multiple Sclerosis.

The model proposed in \cite{lombardo2017demyelination} is a PDE system describing the interaction between the density of activated immune cells (macrophages) $m = m(t, x)$, the concentration of chemical species (cytokine)  $c = c(t, x)$ secreted by the immune cells, and the density of the destroyed oligodendrocytes $d = d(t, x)$. Proper Turing instability analysis of the proposed reaction--diffusion system is performed, showing the formation of spatial patterns when the chemotactic coefficient overcomes a proper threshold; two-dimensional numerical simulations show the appearance of different kinds of patterns (as circular rings, or small clusters, etc.), similar to those really observed in a damaged brain (like Balo's plaques).
Spatial modulation of the Turing--type structures through Eckhaus and zigzag instability is also addressed in \cite{bilotta2018eckhaus}.
Rigorous results on this system have been proven in \cite{desvillettes2020global, Desvillettes-Giunta-RicMat2021}, concerning the global existence of weak and strong solutions, uniform bounds in time for such solutions, and nonlinear stability results under proper assumptions on the chemotaxis term. On the other hand, the stability analysis fits the setting of PDEs--ODEs systems of the recent paper \cite{kowall2022stability}. 

In this paper, we consider a modification of the model proposed in \cite{lombardo2017demyelination}, changing the term describing the production and saturation of the activated macrophages. In \cite{lombardo2017demyelination}, the rate of macrophages activation is modelled by means of a logistic functional form describing the proliferation and saturation effects taken into account also in previous ODE models for acute inflammation \cite{Kumar-etal-2004, Reynolds-etal-2006}.
However, since it is commonly conjectured that activation of microglia is responsible for the appearance of Multiple Sclerosis lesions \cite{Ponomarev-etal-2005}, the underlying mechanism still remains unknown \cite{lombardo2017demyelination}. It is also plausible that a minimal amount of activated macrophages is needed to start and sustain the inflammatory response. This case was not considered in \cite{lombardo2017demyelination}, since they modelled the macrophage activation using a logistic growth function.
In the present work, we thus take into account a different model for macrophages activation, by including the so-called Allee effect \cite{allee-oikos, allee1949principles}, which is a growth function used in population dynamics to take into account undercrowding effects (see \cite{buffoni2011effects, buffoni2016dynamics, rebelo2019persistence} and references therein).
Mainly two types of Allee effects have been considered in the literature \cite{Hastings}, namely the weak and the strong (or critical) Allee effect, depending on the fact that the prey growth function is non-negative or negative, respectively, for small population densities. Examples of the use of the Allee effect in biomedical applications can be found for instance in  \cite{reisch2019chemotactic}, to model in a liver inflammation the fact that a very small virus population can be eliminated
by local immune reactions, and not necessarily increases as prescribed by a logistic growth.

The aim of this work is to investigate how the choice of a different model for macrophage activation can change the pattern scenarios presented and discussed in \cite{lombardo2017demyelination}. To this end, we will perform a weakly nonlinear analysis of the modified reaction--diffusion system, together with a bifurcation analysis.
The weakly nonlinear approach has been adopted in \cite{VanHecke-etal, Gambino-etal-2012, barresi2016wavefront, lombardo2017demyelination} and it is based on the multiple scales method; suitable expansions in terms of a parameter measuring the dimensionless distance with respect to a suitable critical bifurcation value allow to derive amplitude equations which yield the form of the pattern close to the bifurcation threshold. We perform this asymptotic procedure up to the third-order accuracy, getting a Stuart–Landau equation for the amplitude of spatially periodic solutions. 
On the other hand, we also perform a bifurcation analysis by computing the bifurcation structure and its deformation under parameter variations. In this regard, we exploit the numerical continuation software \texttt{pde2path} \cite{uecker2014pde2path, rademacher2018oopde, uecker2021numerical}, based on a FEM discretisation of the stationary problem. The software has been recently used beyond its standard setting, for instance in cross-diffusion \cite{CKCS, breden2021influence, soresina2022hopf} and fractional diffusion problems \cite{ehstand2021numerical}.

The paper is organised as follows. In Section \ref{sec:MathMod}, the mathematical model with the Allee effect is introduced; in Sections \ref{sec:TI} and \ref{sec:wna}, we perform the Turing instability and the weakly nonlinear analysis, respectively. Numerical results, namely bifurcation diagrams obtained via \texttt{pde2path} and numerical simulations, are shown in Section \ref{sec:NumSim}. Finally, some concluding remarks can be found in Section \ref{sec:conclusion}. The Matlab scripts for the numerical continuation are freely accessible in the GitHub folder \cite{GitHubFolder}.

\section{The mathematical model}\label{sec:MathMod}
We want to describe the initial stage of the disease, which is characterized by the production of pro-inflammatory cytokines $\tilde{c}$ by macrophages and activated microglia $\tilde{m}$ and by the emergence of lesions due to apoptosis of oligodendrocyte $\tilde{d}$. We suppose that all the concentrations depend on the position $\mathbf{X}\in\Omega\subset\mathbb{R}^n$ and time $T\in\mathbb{R}^+$.

We consider the following system of partial differential equations
\begin{equation*}
\begin{aligned}
&\dfrac{\partial \tilde{m}}{\partial T}\,=\,D\Delta_\mathbf{X}\tilde{m}+\tilde{G}_j(\tilde{m})-\nabla_\mathbf{X}\cdot(\Psi(\tilde{m})\nabla_\mathbf{X}\tilde{c})\,, &\textnormal{on}\,\Omega\times\mathbb{R}^+\,,\\
&\dfrac{\partial \tilde{c}}{\partial T}\,=\,\dfrac{1}{\nu}\,[\varepsilon\Delta_\mathbf{X}\tilde{c}+\mu\tilde{d}-\alpha\tilde{c}+b\tilde{m}]\,,&\textnormal{on}\,\Omega\times\mathbb{R}^+\,,\\
&\dfrac{\partial \tilde{d}}{\partial T}\,=\,\kappa\tilde{F}(\tilde{m})\tilde{m}(\bar{d}-\tilde{d})\,,&\textnormal{on}\,\Omega\times\mathbb{R}^+\,,\\
&\dfrac{\partial \tilde{m}}{\partial n}=\dfrac{\partial \tilde{c}}{\partial n}=\dfrac{\partial \tilde{d}}{\partial n}=0\,, &\textnormal{on}\,\partial\Omega\times\mathbb{R}^+\,,\\
&\tilde{m}(\mathbf{X},0)=\tilde{m}_{in}(\mathbf{X}),\, \tilde{c}(\mathbf{X},0)=\tilde{c}_{in}(\mathbf{X}),\, \tilde{d}(\mathbf{X},0)=\tilde{d}_{in}(\mathbf{X}),\,&\textnormal{on}\,\Omega\,,
\end{aligned}
\end{equation*}
being $\Omega$ a bounded and connected domain. The evolution of macrophages/microglia density is governed by random movements with constant diffusion coefficient $D$, a production/decay term $\tilde{G}_j$ (where $j=1,2$ denotes two different functions detailed below) and a chemotactic motion towards regions with a high concentration of chemoattractants; the function $\Psi(\tilde{m})=\psi\tilde{m}/(\bar{m}+\tilde{m})$, where $\bar{m}$ is the characteristic density of macrophages and $\psi$ is the maximal chemotactic rate. The cytokines evolution is ruled by random diffusion with diffusion coefficient $\varepsilon$ and by production/decay terms with constant proliferation coefficients $\mu,\, b$ and death coefficient $\alpha$. For what concerns oligodendrocytes, their evolution is governed by a production term, where $\kappa$ measures the destructive strength of the macrophages, $\tilde{F}(\tilde{m})=\tilde{m}/(\bar{m}+\tilde{m})$ and $\bar{d}$ is the characteristic density of oligodendrocytes.\\
Regarding the production term for macrophages $\tilde{G}_j$, we shall consider, contrast and compare the effects of two different growth functions; the first choice is a logistic term
\begin{equation*}
	\tilde{G}_1(\tilde{m})=\lambda\tilde{m}(\bar{m}-\tilde{m})\,,
\end{equation*}
where $\lambda$ is the production rate; alternatively, the production rate can also depend on the concentration of macrophages and we shall consider a cubic function
\begin{equation*}
	\tilde{G}_2(\tilde{m})=\tilde{\lambda}\tilde{m}(\bar{m}-\tilde{m})(\tilde{m}-\hat{m})\,,
\end{equation*}
where $\tilde{\lambda}$ is the rate of increase in presence of an Allee effect, described by the additional parameter $\hat{m}<\bar{m}$.  When $\hat{m}\leq 0$, the Allee effect is weak, namely the growth rate is reduced but still positive; when $\hat{m}> 0$, we are in presence of a strong Allee effect, namely $\hat{m}$ describes a threshold in the concentration of macrophages allowing their growth.

We introduce the scaled variables and parameters
\begin{equation*}
	\begin{aligned}
		m=\tilde{m}/\bar{m}\;,\;d = \tilde{d}/\bar{d}\;,\;c&=\dfrac{\alpha}{b\bar{m}}\,\tilde{c}\;,\;t=\lambda\bar{m}T\;,\;\mathbf{x}=\sqrt{\dfrac{\lambda\bar{m}}{D}}\,\mathbf{X}\,,\\
		\chi=\dfrac{\bar{b}}{\alpha D}\,\psi\;,\;\tau=\dfrac{\lambda \bar{m}}{\alpha}\,\nu\;,\;\epsilon&=\dfrac{\lambda\bar{m}}{\alpha D}\,\varepsilon\;,\;\beta=b/\bar{b}\;,\;r=\kappa/\lambda\;,\;\delta=\dfrac{\bar{d}}{\bar{m}\bar{b}}\,\mu\,,
	\end{aligned}
\end{equation*}
where $\bar{b}$ is a typical production rate for macrophages.\\
The governing equations can be rewritten in the nondimensional form 
\begin{equation}
\begin{aligned}
&\dfrac{\partial m}{\partial t}\,=\,\Delta_\mathbf{x}m+G_j(m)-\nabla_\mathbf{x}\cdot(\Phi(m)\nabla_\mathbf{x}c)\,,&\textnormal{on}\,\Omega\times\mathbb{R}^+\,,\\
&\dfrac{\partial c}{\partial t}\,=\,\dfrac{1}{\tau}\,[\epsilon\Delta_\mathbf{x}c+(\delta d-c+\beta m)]\,,&\textnormal{on}\,\Omega\times\mathbb{R}^+\,,\\
&\dfrac{\partial d}{\partial t}\,=\,r F(m)m(1-d)\,,&\textnormal{on}\,\Omega\times\mathbb{R}^+\,,\\
&\dfrac{\partial m}{\partial n}=\dfrac{\partial c}{\partial n}=\dfrac{\partial d}{\partial n}=0\,, &\textnormal{on}\,\partial\Omega\times\mathbb{R}^+\,,\\
&m(\mathbf{x},0)=m_{in}(\mathbf{x}),\, c(\mathbf{x},0)=c_{in}(\mathbf{x}),\, d(\mathbf{x},0)=d_{in}(\mathbf{x}),\,&\textnormal{on}\,\Omega\,,
\end{aligned}
\label{nondim_syst}
\end{equation}
with $\chi,\tau,\epsilon,r>0$ and $\beta,\delta\ge 0$; the production term $G_j,\, j=1,2$ for logistic and cubic growth function has the following form
\begin{equation*}
	G_1(m)=m(1-m)\;\text{ and }\;G_2(m)=\Lambda m(1-m)(m-M)\,,
\end{equation*}
respectively, where $\Lambda$ is a proper nondimensional production rate and $M=\hat{m}/\bar{m}<1$.

The rigorous results proven in \cite{desvillettes2020global, Desvillettes-Giunta-RicMat2021} for the model in \cite{lombardo2017demyelination} concerning the global existence of weak and strong solutions, uniform bounds in time for such solutions, also apply in this case. 

\section{Turing instability}\label{sec:TI}
In this section, we analyse the linear stability of system \eqref{nondim_syst} by focusing on the formation of stationary patterns due to the destabilisation of a homogeneous steady state. First, we look for equilibria of the homogeneous system
\begin{equation*}
	\begin{aligned}
		\dfrac{\partial m}{\partial t}\,&=\,G_j(m)\,,\\
		\dfrac{\partial c}{\partial t}\,&=\,\dfrac{1}{\tau}(\delta d-c+\beta m)\,,\\
		\dfrac{\partial d}{\partial t}\,&=\,r F(m)m(1-d)\,.
	\end{aligned}
\end{equation*}
We observe that the line $\{(0,\delta d,d), \,d\in\mathbb{R}^+\}$ of steady states characterized by the absence of macrophages is stable for both logistic growth and a weak Allee effect, and unstable for a strong Allee effect. The coexistence state $P_*=(m_*,c_*,d_*)=(1,\beta+\delta,1)$ is a stable equilibrium in all cases (logistic growth, weak and strong Allee effect). When a strong Allee effect is considered, we have an additional unstable steady state $P_\sharp=(m_\sharp,c_\sharp,d_\sharp)=(M,\beta M+\delta,1)$.\\
Now we investigate the conditions leading to the destabilisation of $P_*$ in presence of spatial diffusion  and the chemotactic term.\\
We denote by $\boldsymbol{\bar{w}}=(w_m,w_c,w_d)=(m-m_*,c-c_*,d-d_*)$ the perturbation with respect to the homogeneous state and we consider for it the following linearised problem
\begin{equation*}
	\dfrac{\partial \boldsymbol{w}}{\partial t}\,=\,\mathbf{J}^\prime\boldsymbol{\bar{w}}+\mathbf{D}^\prime\Delta_\mathbf{x}\boldsymbol{\bar{w}}\,,
\end{equation*}
where, for $j=1,2$ we have
\begin{equation}
	\mathbf{J}^\prime\,=\,
	\begin{pmatrix}
		-[\Lambda(1-M)]^{j-1} & 0 & 0\\
		\\
		\dfrac{\beta}{\tau} & -\dfrac{1}{\tau} & \dfrac{\delta}{\tau}\\
		\\
		0 & 0 & -\dfrac{r}{2}
	\end{pmatrix}
	\,,\quad
	\mathbf{D}^\prime\,=\,
	\begin{pmatrix}
		1 & -\dfrac{\chi}{2} & 0\\
		\\
		0 & \dfrac{\epsilon}{\tau} & 0\\
		\\
		0 & 0 & 0
	\end{pmatrix}\,.
	\label{JD}
\end{equation}
From now on, we assume that $\Omega$ is a rectangular domain in $\mathbb{R}^n$; we look for solutions of the form $\boldsymbol{\bar{w}}\,\propto \text{exp}(\lambda t+i \mathbf{k}\cdot\mathbf{x})$  and compute the eigenvalues $\lambda$. One is easily given by $-r/2$; the others are given as proper functions of $k=\lvert\mathbf{k}\rvert$, by solving the dispersion relation 
\begin{equation*}
	\det(\lambda\mathbf{I}-\mathbf{J}+k^2\mathbf{D})\,=\,0\,,
\end{equation*}
where
\begin{equation}
	\mathbf{J}=
	\begin{pmatrix}
		-[\Lambda(1-M)]^{j-1} & 0\\
		\\
		\dfrac{\beta}{\tau} & -\dfrac{1}{\tau}	
	\end{pmatrix}\,\text{ and }\,
	\mathbf{D}=\mathbf{D}(\chi)=
	\begin{pmatrix}
		1 & -\dfrac{\chi}{2}\\
		\\
		0 & \dfrac{\epsilon}{\tau}	
	\end{pmatrix}\,.
	\label{JD_red}
\end{equation}
We obtain
\begin{equation*}
	\lambda^2+g(k^2)\lambda+h(k^2)\,=\,0\,,
\end{equation*}
where
\begin{equation*}
	g(k^2)\,=\,k^2\text{tr}(\mathbf{D})-\text{tr}(\mathbf{J})\,,\quad h(k^2)\,=\,\det(\mathbf{D})k^4+qk^2+\det(\mathbf{J})\,,
\end{equation*}
being
\begin{equation*} q\,=\,\dfrac{2(1+\xi^2)-\chi\beta}{2\tau}\,, \quad \xi=\sqrt{[\Lambda(1-M)]^{j-1}\epsilon}\,.
\end{equation*}
The Turing instability is guaranteed by a positive eigenvalue  for at least one $k>0$. Since for all the values of the parameters
\begin{equation*}
	g(k^2)\,=\,k^2(1+\epsilon/\tau)+[\Lambda(1-M)]^{j-1}+1/\tau\,>\,0\,,
\end{equation*}
the presence of a positive eigenvalue is guaranteed by $h<0$ for some $k$. 
Being $h$ a second-order polynomial in $k^2$ with positive highest and lowest order coefficients, we require the minimum of $h$ to be negative. The minimum is reached at
\begin{equation}
	k^2_c\,=\,-\dfrac{q}{2\det(\mathbf{D})}\,=\,-\dfrac{2(1+\xi^2)-\chi\beta}{4\epsilon}\,.
	\label{eq:k2formula}
\end{equation}
For this to be positive, it is required
\begin{equation*}
	2(1+\xi^2)-\chi\beta\,<\,0\,.
\end{equation*}
Here we see that it is not diffusion but the chemotactic term the key ingredient destabilising the homogeneous steady state. We consider $\chi$ as the bifurcation parameter, obtaining the threshold value
\begin{equation*}
	\chi\,>\,\bar{\chi}\,=\,\dfrac{2(1+\xi^2)}{\beta}\,.
\end{equation*}
We look for a critical value for this parameter by imposing
\begin{equation*}
	h(k_c^2)\,=\,0\,,
\end{equation*}
and we find the threshold for the bifurcation parameter $\chi$
\begin{equation}
	\chi_c\,=\,\dfrac{2(1+\xi)^2}{\beta}\,>\,\bar{\chi}
	\label{chi_c}
\end{equation}
and hence, substituting $\chi_c$ in \eqref{eq:k2formula}, we have
\begin{equation}
	k_c^2\,=\,\dfrac{\xi}{\epsilon}=\sqrt{\dfrac{[\Lambda(1-M)]^{j-1}}{\epsilon}}\,.
	\label{k_c}
\end{equation}
Note that $\chi_c$ decreases as $M$ or $\beta$ increases, and it increases as $\epsilon$ or $\Lambda$ increases.

Bifurcation points may occur for $\chi>\chi_c$. More precisely, the mode associated to the eigenvalues $k^2$ of the Laplacian destabilises at $\chi=\chi_{A,k^2}$, where
\begin{equation}
\chi_{A,k^2}=\dfrac{\epsilon k^4+\tau\left(2+\epsilon[\Lambda(1-M)]^{j-1}\right)k^2+[\Lambda(1-M)]^{j-1}}{\beta k^2},
\label{eq:chi_bif}
\end{equation}
obtained by solving $h(k^2)=0$ with respect to $\chi$. Note that, in the logistic case, this reduces to
\begin{equation*}
\chi_{L,k^2}=\dfrac{\epsilon k^4+\tau\left(2+\epsilon\right)k^2+1}{\beta k^2}.
\end{equation*}

\section{Weakly nonlinear analysis}\label{sec:wna}

We are now interested in deriving the Stuart--Landau \cite{lombardo2017demyelination, VanHecke-etal} equation for the amplitude of spatially periodic solutions of system \eqref{nondim_syst} in the one-dimensional case ($n=1$) in the domain $\Omega=[0,\bar{L}]$. We use a multiple-scale method for performing a weakly nonlinear analysis around the stationary state $P_*$. Since there is no diffusion in the equation for $d$, we shall consider the reduced system in the variables $\boldsymbol{w}=(w_m,w_c)$. Anyway, repeating the same procedure as follows for the whole set of variables $\boldsymbol{w}=(w_m,w_c,w_d)$, it is easy to show that this gives a null solution for $w_d$. 
We rewrite the reduced system as
\begin{equation}
	\dfrac{\partial \boldsymbol{w}}{\partial t}\,=\,\mathcal{L}_\chi\boldsymbol{w}+\mathcal{N}\boldsymbol{w}\,,
	\label{lin_syst}
\end{equation}
where the linear operator is defined as $\mathcal{L}_\chi=\mathbf{J}+\mathbf{D}\partial_{xx}$ ($\mathbf{J}$ and $\mathbf{D}$ are given in \eqref{JD_red}) while $\mathcal{N}$ keeps track of nonlinear contributions.

\begin{figure}
\centerline{\includegraphics[width=0.7\textwidth]{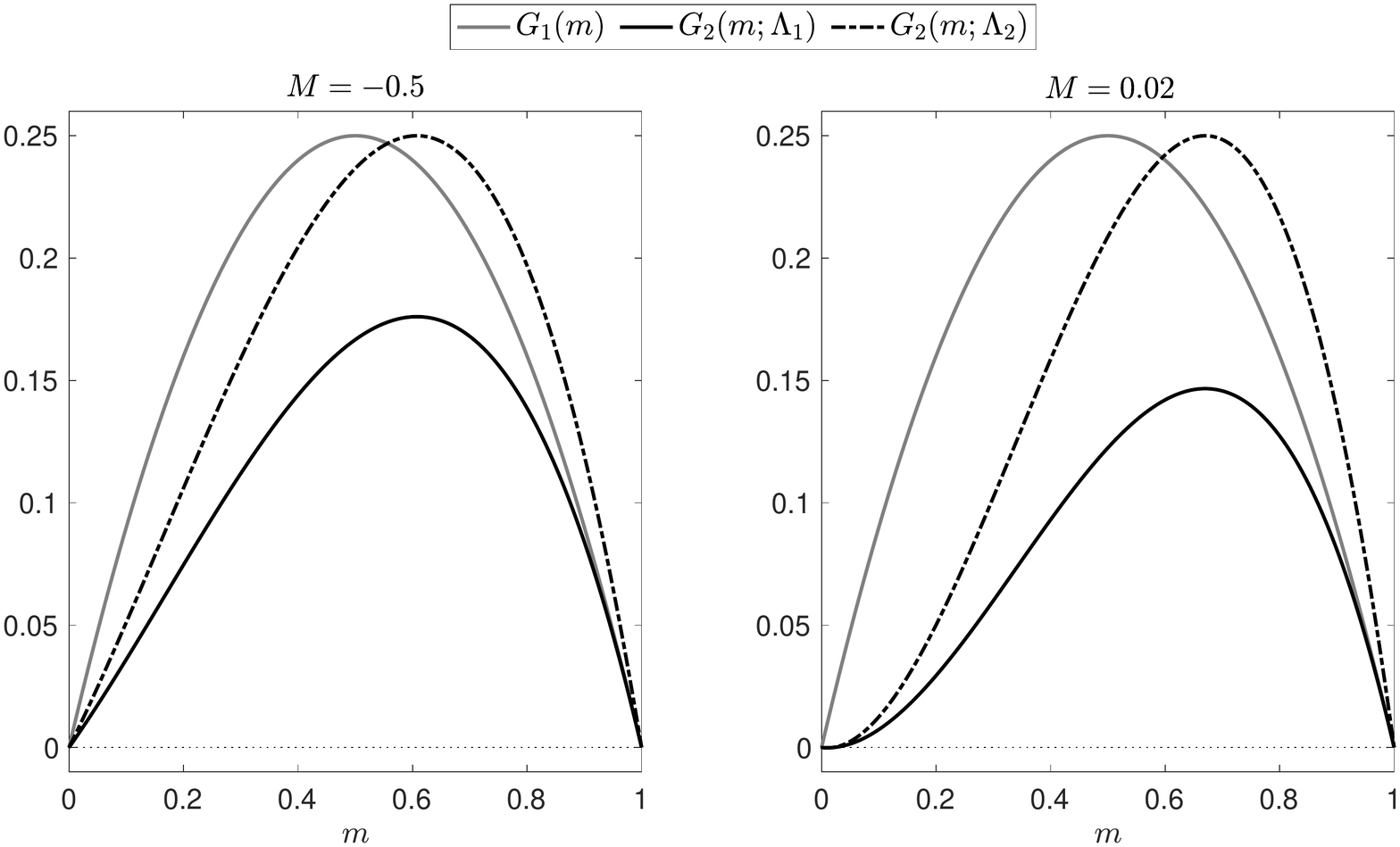}}
\caption{Comparison of the profile of logistic growth function $G_1$ (grey line) with cubic growth functions (black lines) with parameter $\Lambda$ as in Case 1 ($G_2(m;\Lambda_1)$, black solid line) and in Case 2 ($G_2(m;\Lambda_2)$, dashed-dotted line). The left panel corresponds to the weak Allee effect ($M=-0.5$), while the right panel to a strong Allee effect ($M=0.02$).}
\label{fig:growth_cases1&2}
\end{figure}

We introduce a small parameter $\eta^2=(\chi-\chi_c)/\chi_c$ measuring the distance of the bifurcation parameter from its critical value; we consider the following expansions
\begin{equation*}
	\begin{aligned}
		\boldsymbol{w}&=\eta\boldsymbol{w}_1+\eta^2\boldsymbol{w}_2+\eta^3\boldsymbol{w}_3+O(\eta^4)\,,\\
		\chi&=\chi_c+\eta^2\chi_2+O(\eta^4)\,,
	\end{aligned}
\end{equation*}
where $\boldsymbol{w}_i=(w_m^{(i)},w_c^{(i)}), \, i=1,2,3$; for what concerns the time dependence of the solutions, we assume a multiple scale dependence $\boldsymbol{w}=\boldsymbol{w}(T_2,T_4,\ldots)$, where $T_{2i}=\eta^{2i} t, \, i=1,2$, so that the time derivative can be expanded as follows
\begin{equation*}
	\dfrac{\partial}{\partial t}=\eta^2\dfrac{\partial}{\partial T_2}+\eta^4\dfrac{\partial}{\partial T_4}+O(\eta^5)\,.
\end{equation*}
We can substitute the previous expansions in (\ref{lin_syst}) and collect the terms of the same order.\\
At the first order of accuracy, we have to solve
\begin{equation}
	\mathcal{L}_{\chi_c}\boldsymbol{w}_1=\boldsymbol{0}\,;
	\label{syst_ord_1}
\end{equation}
by imposing no flux boundary conditions, we look for solutions of the form
\begin{equation}
	w_\ell^{(1)}\,=\,\rho_\ell A(T_2,\ldots)\cos(k_cx)\,,\qquad\ell=m,c\,.
	\label{forma_ord_1}
\end{equation}
System \eqref{syst_ord_1} admits non trivial solutions; in particular, the following condition holds
\begin{equation*}
\rho_m=\dfrac{1+\xi}{\beta}\,\rho_c\,.
\end{equation*}
At the second order of accuracy, we have to solve a non-homogeneous system
\begin{equation}
	\mathcal{L}_{\chi_c}\boldsymbol{w}_2=\boldsymbol{F}(\boldsymbol{w}_1)\,,
	\label{syst_ord_2}
\end{equation}
where $\boldsymbol{F}$ contains only constants and terms proportional to $\cos (2 k_c x)$.
For the Fredholm alternative, the solvability condition for \eqref{syst_ord_2} is $<\boldsymbol{F},\boldsymbol{\psi}>=0$, where $\boldsymbol{\psi}\in\ker(\mathcal{L}^*)$ and $\mathcal{L}^*$ is the adjoint operator of $\mathcal{L}_{\chi_c}$, and it is automatically satisfied, since $\boldsymbol{\psi}=(\psi_1\cos(k_cx),\psi_2\cos(k_cx))$ with
\begin{equation*}
	\psi_2\,=\,\dfrac{\tau}{\beta\epsilon}\,\xi(1+\xi)\psi_1\,.
\end{equation*}
The solution can be easily written as
\begin{equation}
	w^{(2)}_\ell\,=\,A^2\left[\mu_\ell+\theta_\ell\cos(2k_cx)\right]\,,\qquad \ell=m,c\,,
	\label{momenti_ordine_2}
\end{equation}
with
\begin{equation*}
	\begin{aligned}
		\mu_m\,&=\,-\dfrac{\rho_m^2}{2}\,\left(1+\omega\right),&\qquad\theta_m\,&=\,-\dfrac{\rho_m^2}{18}\left(\omega-\dfrac{1}{\xi}\right)(1+4\xi)&\\[0.2cm]
		\mu_c\,&=\,-\beta\,\dfrac{\rho_m^2}{2}\,\left(1+\omega\right),&\qquad\theta_c\,&=\,-\beta\,\dfrac{\rho_m^2}{18}\left(\omega-\dfrac{1}{\xi}\right)&
	\end{aligned}
\end{equation*}
and $\omega=(j-1)/(1-M)$.
At the third order, we have to solve
\begin{equation}
	\mathcal{L}_{\chi_c}\boldsymbol{w}_3=\boldsymbol{G}(\boldsymbol{w}_1,\boldsymbol{w}_2)\,,
	\label{syst_ord_3}
\end{equation}
where the right hand side can be rewritten in the form
\begin{equation}
	\boldsymbol{G}(\boldsymbol{w}_1,\boldsymbol{w}_2)=\left[\mathbf{G}^0\dfrac{dA}{dT_2}+\mathbf{G}^1A+\mathbf{G}^2A^3\right]\cos(k_cx)+\mathbf{G}^*\,,
	\label{RHS_ordine_3}
\end{equation}
with
\begin{equation*}
	\begin{aligned}
		\mathbf{G}^0=(G_m^0,G_c^0)&=\rho_m\left(1,\dfrac{\beta}{1+\xi}\right)\\
		\mathbf{G}^1=(G_m^1,G_c^1)&=\rho_m\left(-\dfrac{\chi_2}{2}\dfrac{\beta}{\epsilon}\dfrac{\xi}{1+\xi},0\right)\\
		\mathbf{G}^2=(G_m^2,G_c^2)&=\Bigg(2\dfrac{\xi^2}{\epsilon}\left(1+\omega\right)\rho_m\left(\mu_m+\dfrac{\theta_m}{2}\right)+\dfrac{3}{4\epsilon}\omega\xi^2\rho_m^3\\
		&\hspace{1cm}-\dfrac14\chi_ck_c^2\left[\rho_m\theta_c+\rho_c\left(\mu_m-\dfrac{\theta_m}{2}\right)-\dfrac18\rho_m^2\rho_c\right],0\Bigg)
	\end{aligned}
\end{equation*}
and $\mathbf{G}^*\propto\cos(3k_cx)$ satisfies the Fredholm condition. By imposing the solvability condition at this order, we obtain the Stuart--Landau equation for the amplitude
\begin{equation*}
	\dfrac{dA}{dT_2}\,=\,\sigma A-LA^3\,,
\end{equation*}
where
\begin{equation*}
	\sigma\,=\,-\dfrac{G^1_m\psi_1+G^1_c\psi_2}{G^0_m\psi_1+G^0_c\psi_2}\,\text{ and }\,L\,=\,\dfrac{G^2_m\psi_1+G^2_c\psi_2}{G^0_m\psi_1+G^0_c\psi_2}\,.
\end{equation*}
We can observe that 
\begin{equation*}
	\sigma=\dfrac12\beta\chi_2\dfrac{\xi}{(1+\xi)(\epsilon+\tau\xi)}>0\,;
\end{equation*}
in fact, we consider the perturbation $\chi_2>0$ because Turing instability may occur only for $\chi>\chi_c$. 
For what concerns $L$, it can be rewritten as
\begin{equation*}
	L=-\dfrac{\rho_m^2}{144(\epsilon+\tau\xi)}p(\xi,\omega)\,,
\end{equation*}
and its sign depends on the sign of the polynomial
\begin{equation*}
	p(\xi,\omega)=4\omega(8\omega+9)\xi^3+(152\omega^2+122\omega+63)\xi^2-(46\omega+55)\xi+2\,.
\end{equation*}
In absence of Allee effect, i.e.~when we consider the logistic function $G_1(m)$, we have
\begin{equation*}
	\xi=\sqrt{\epsilon}\;,\;\omega=0\,,
\end{equation*}
and the polynomial $p$ reduces to
\begin{equation*}
	p(\sqrt{\epsilon})=63\epsilon-55\sqrt{\epsilon}+2\,;
\end{equation*}
therefore, we obtain that
\begin{equation*}
	L>0\,\Longleftrightarrow\, p< 0\,\Longleftrightarrow\,\epsilon\in(0.0014,0.6972)\,.
\end{equation*}
In presence of Allee effect, namely when we consider the cubic function $G_2(m)=\Lambda m(1-m)(m-M)$, we have that the variables $\xi$ and $\omega$ depend on the choice of parameters $M$ and $\Lambda$. In order to contrast and compare the outcomes of different growth terms $G_j, \,j=1,2$ varying these two parameters, we identified two different ways to set the parameter $\Lambda$. We explore the following cases, whose differences are shown in Figure \ref{fig:growth_cases1&2}:
\begin{itemize}
\item[-] \textbf{Case 1:} for a fixed value of $M<1$, we set the parameter $\Lambda=\Lambda_1:=(1-M)^{-1}$ in order to reproduce the bifurcation values of the logistic case, namely $\chi_{A,k^2}=\chi_{L,k^2}$ (see equation \eqref{eq:chi_bif}).\\

\item[-] \textbf{Case 2:} for a fixed value of $M<1$, we set the parameter $\Lambda$ in order to reproduce the maximal growth rate of the logistic function 
\begin{equation*}
    \max_{0\le m\le 1}G_2(m)=\max_{0\le m\le 1}G_1(m)=\frac{1}{4}\,.
\end{equation*}
This leads to the value
\begin{equation*}
\Lambda_2=\left(4m_2(1-m_2)(m_2-M)\right)^{-1}, \quad m_2=\dfrac{(M+1)+\sqrt{M^2-M+1}}{3} >0 \,.
\end{equation*}
By simple algebra, one gets that for all $M<1$,  the value $m_2 \in (M,1)$, and $m_2 =1$ if and only if $M=1$;  therefore, $\Lambda_2>0$ for all $M<1$. Moreover, taking into account that $4x(1-x)\leq 1 $, it follows that   
\begin{equation*}
	\Lambda_2(1-M)>1,\quad \forall M<1,
\end{equation*}
and hence, for a fixed $k^2$, the bifurcation values $\chi_{A,k^2}$ are greater than the corresponding values for the logistic case, namely $\chi_{A,k^2}>\chi_{L,k^2}$ (see equation \eqref{eq:chi_bif}).
\end{itemize}

\begin{figure}
\centerline{\includegraphics[trim={4cm 2cm 0 0},clip,width=0.7\textwidth]{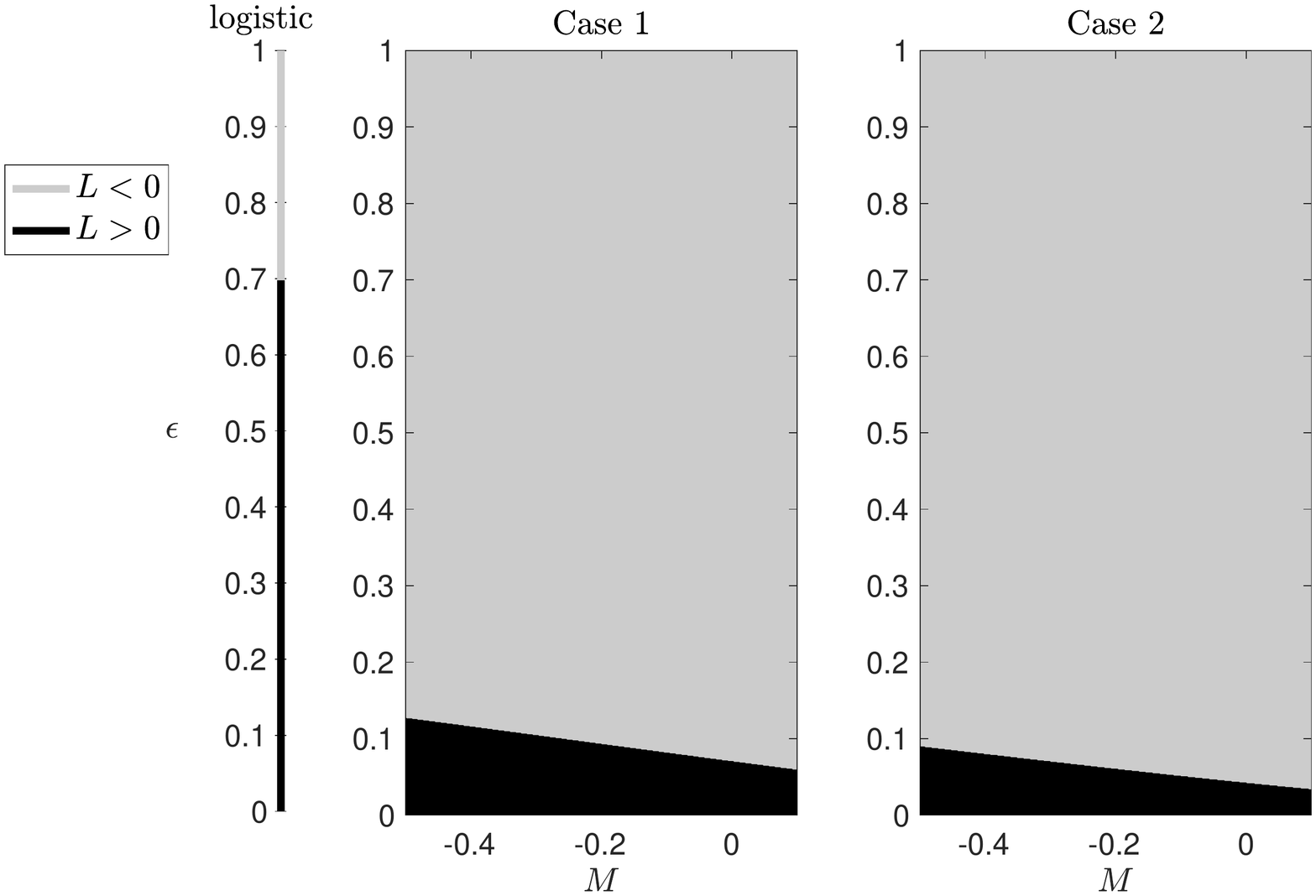}}
\caption{Regions in the $(M,\epsilon)$-plane corresponding to the subcritical case ($L<0$, grey region) and to the supercritical case ($L>0$, black region). The value of $\Lambda$ is set as in Case 1 (left) and as in Case 2 (right). The other parameter values are $\tau=\beta=r=\delta=1$. The thick line on the left corresponds to the logistic case (independent of $M$).}
\label{signL_1e2}
\end{figure}

The regions corresponding to supercritical ($L>0$) and subcritical ($L<0$) bifurcation at $\chi_c$ are reported in Figure \ref{signL_1e2} varying $\epsilon$ for the logistic growth function and in the $(M,\epsilon)$-plane when the Allee effect is considered, in both cases $\Lambda_1$ and $\Lambda_2$.
Both cases show a significant reduction of the range of parameter $\epsilon$ corresponding to the supercritical case ($L>0$), and hence to the existence of a stable equilibrium for the amplitude, especially when we consider a strong Allee effect. Moreover, by confining on realistic values of $\epsilon$, namely $0.5<\epsilon<1.5$ found in \cite{lombardo2017demyelination}, we can observe only a subcritical region in presence of Allee effect; the third order equation cannot capture the amplitude of the pattern and  the analysis should be pushed to higher order ($O(\eta^5)$).

\section{Bifurcation analysis and numerical simulations} \label{sec:NumSim}
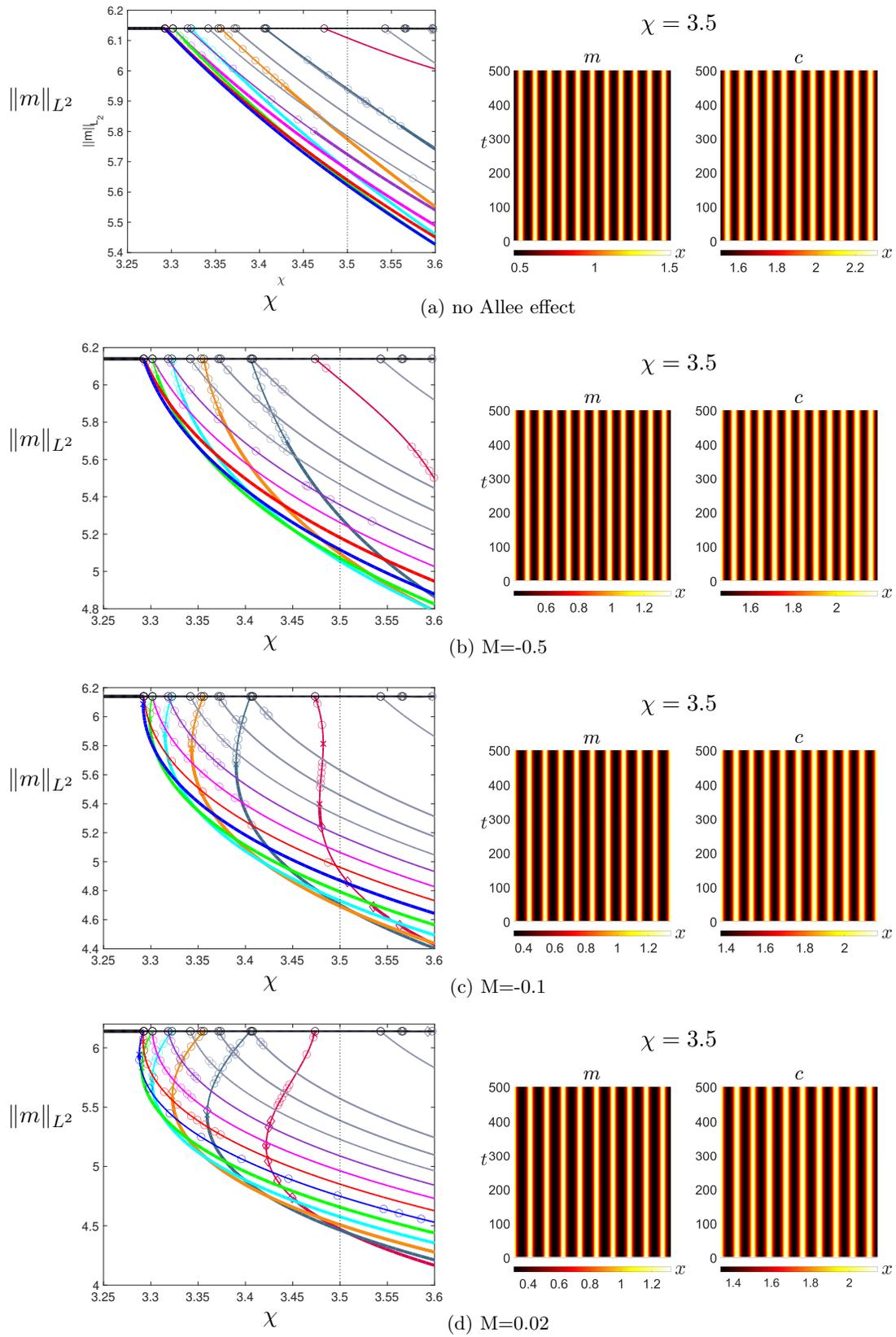
\begin{figure}[!ht]
\centering
\subfloat[no Allee effect
\label{epsilon_0p08_chi_3p50_noallee}]{
\begin{overpic}[width=0.35\textwidth]{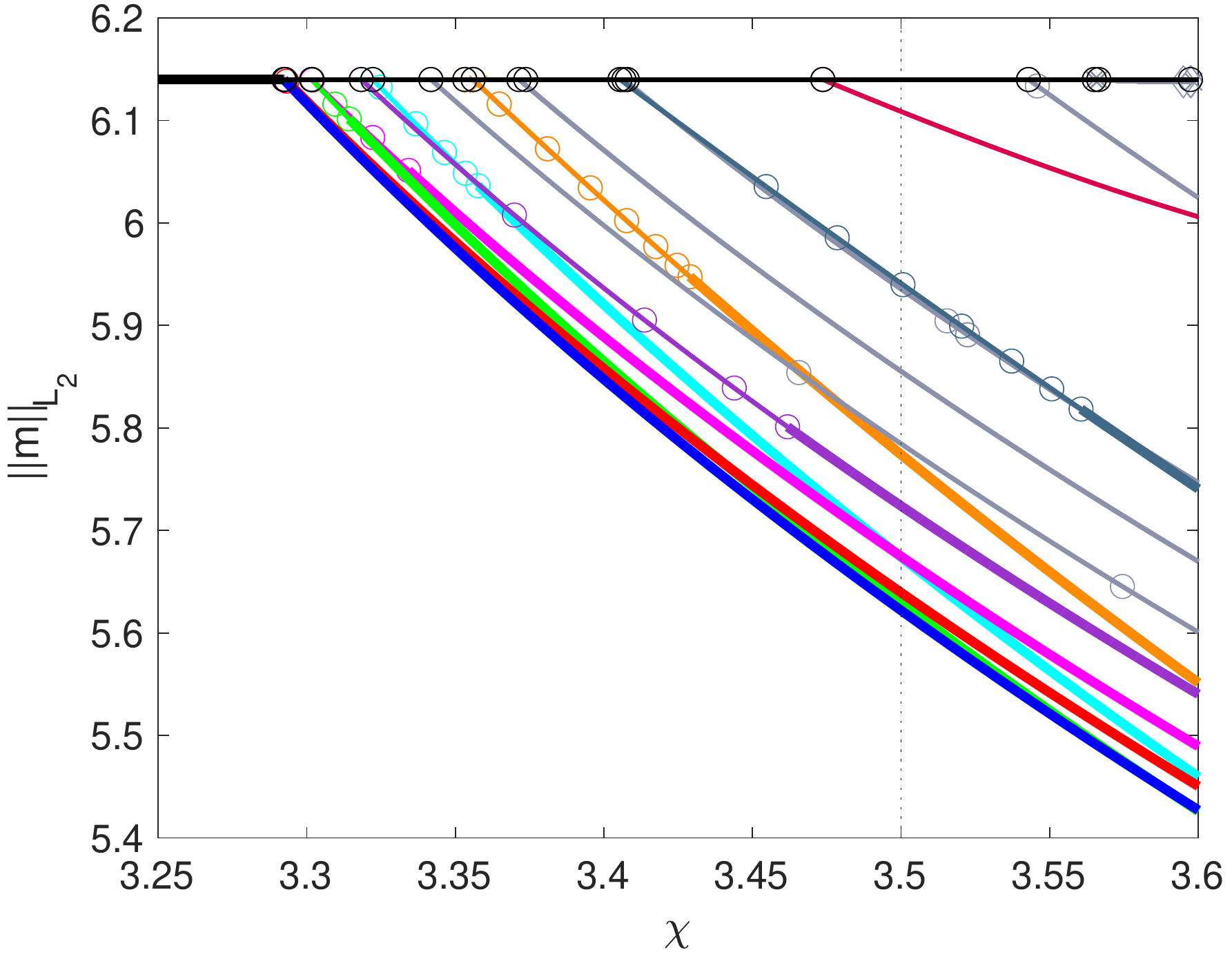}
\put(-20,50){$\|m\|_{L^2}$}
\put(50,-5){$\chi$}
\end{overpic}
\begin{overpic}[width=0.45\textwidth]{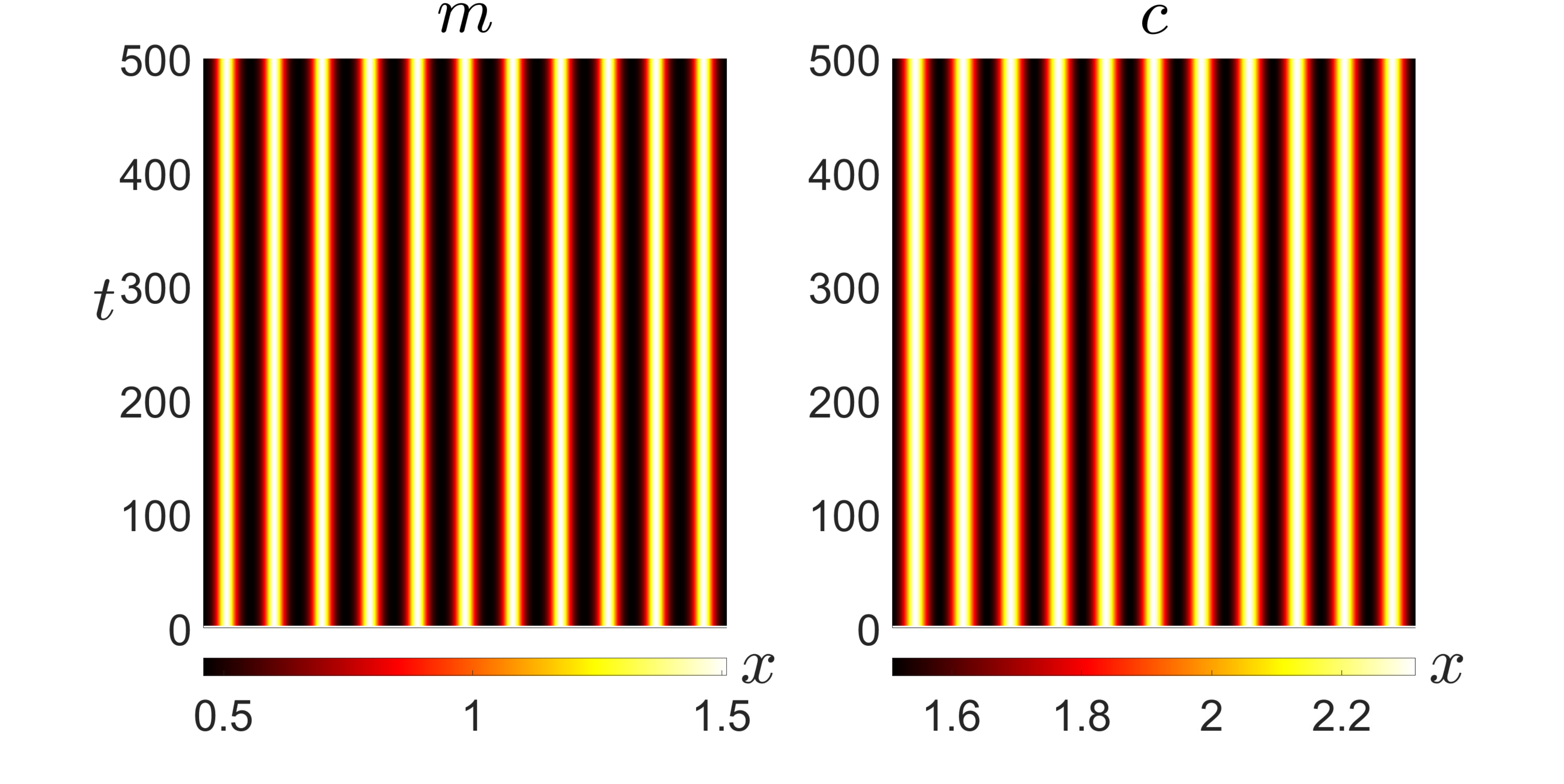}
\put(40,55){$\chi=3.5$}
\end{overpic}
}

\subfloat[M=-0.5
\label{epsilon_0p08_M_m0p5_chi_3p50_caso_2}]{
\begin{overpic}[width=0.35\textwidth]{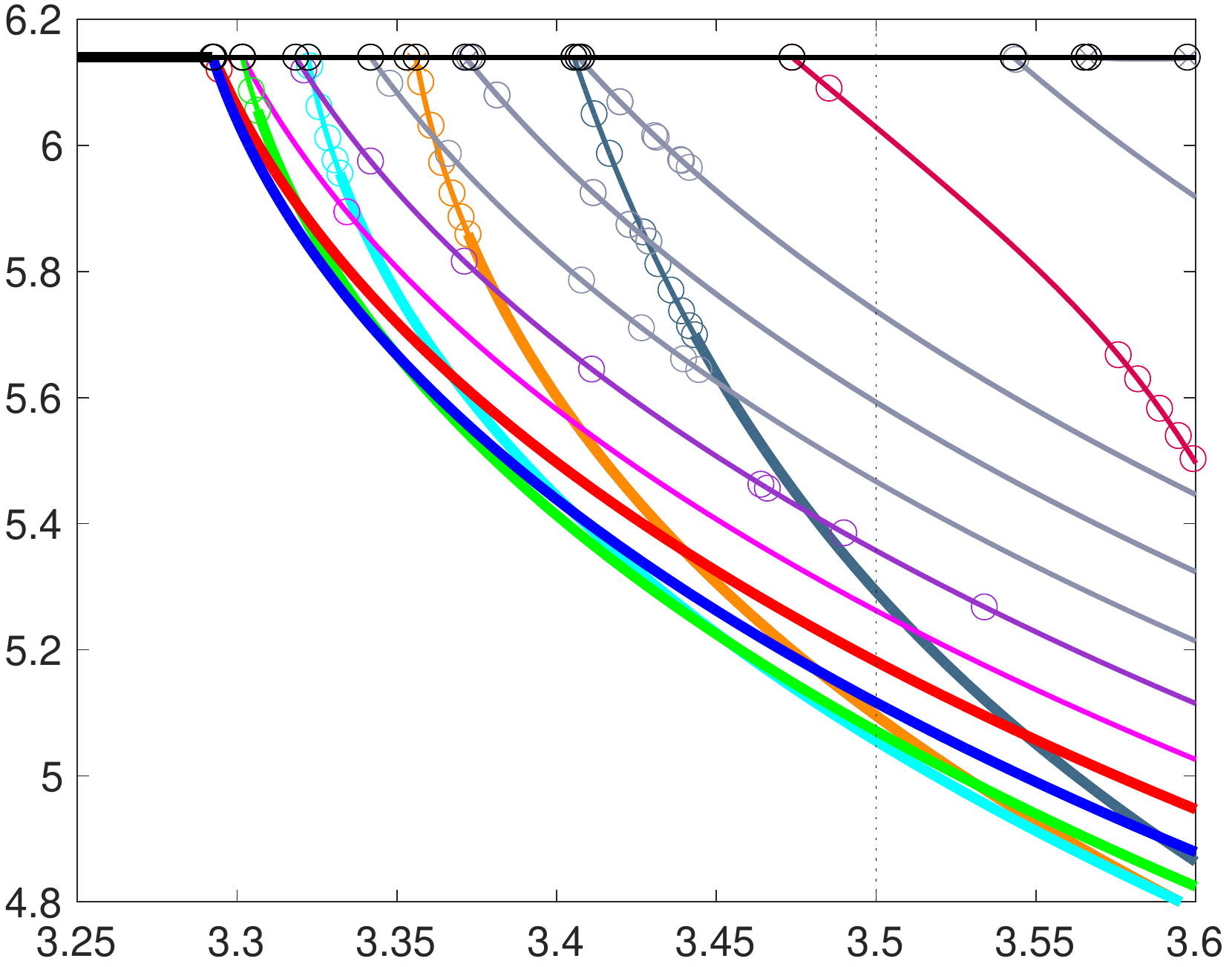}
\put(-20,50){$\|m\|_{L^2}$}
\put(50,-5){$\chi$}
\end{overpic}
\begin{overpic}[width=0.45\textwidth]{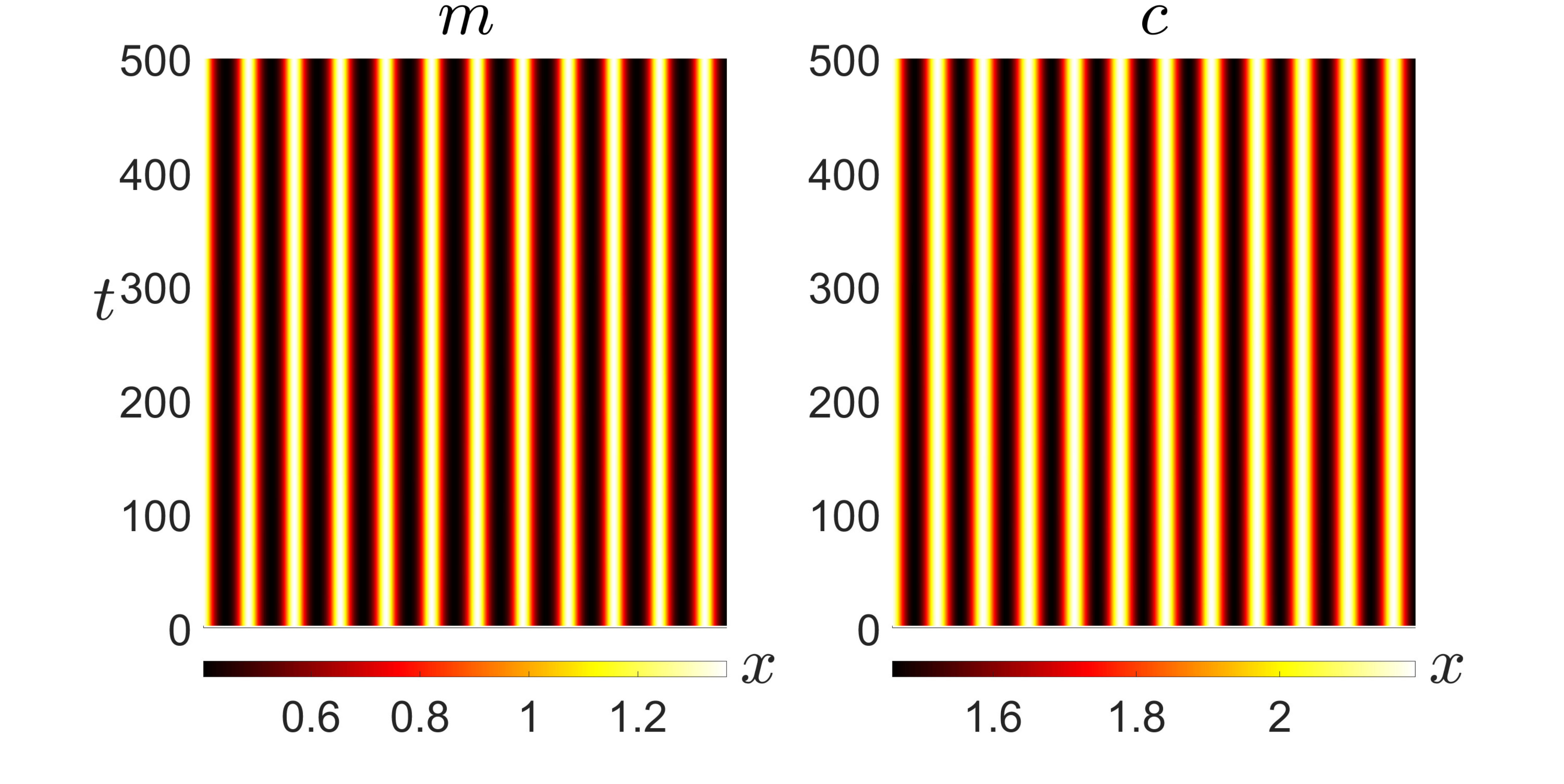}
    \put(40,55){$\chi=3.5$}
\end{overpic}
}

\subfloat[M=-0.1
\label{epsilon_0p08_M_m0p1_chi_3p50_caso_2}]{
\begin{overpic}[width=0.35\textwidth]{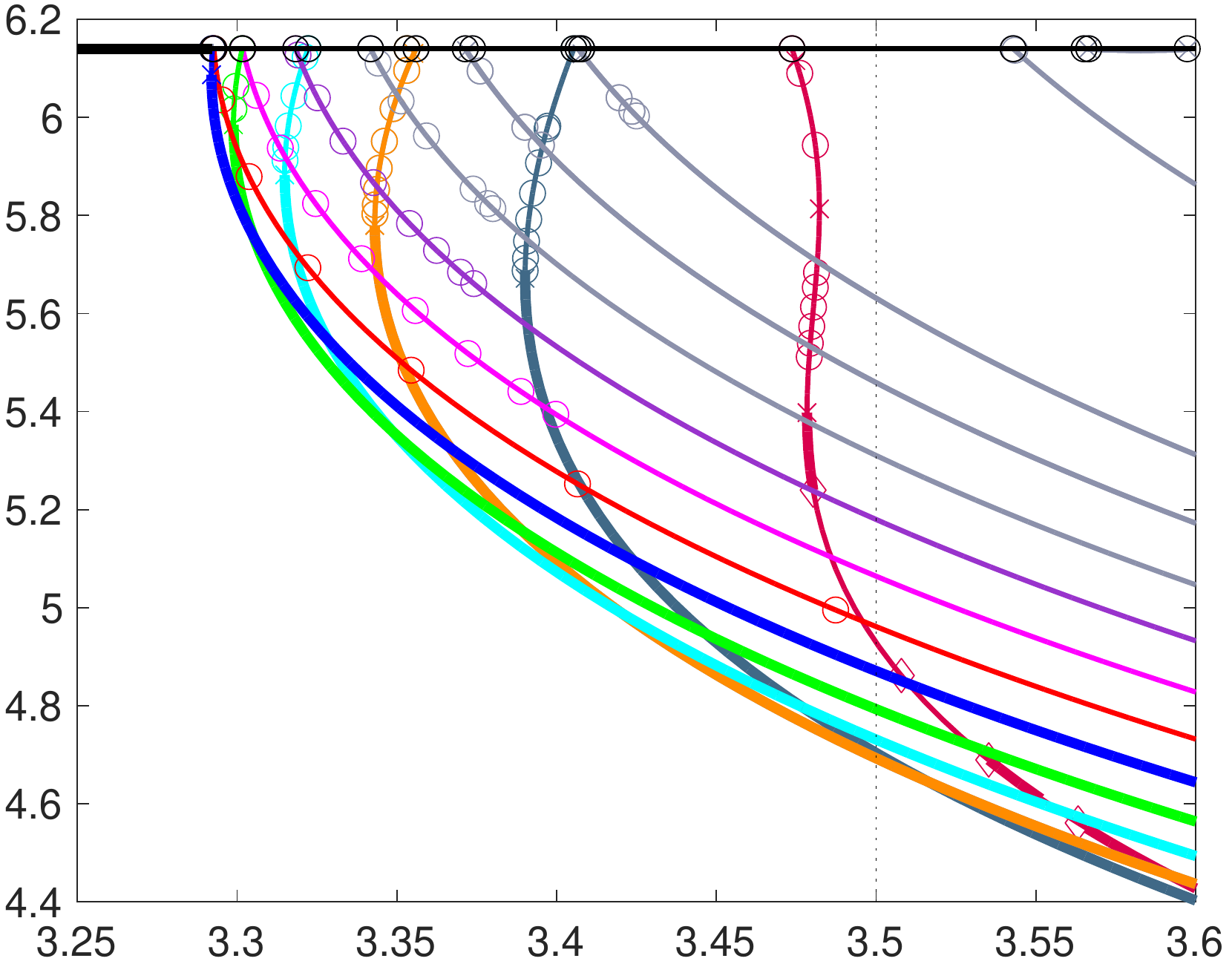}
\put(-20,50){$\|m\|_{L^2}$}
\put(50,-5){$\chi$}
\end{overpic}
\begin{overpic}[width=0.45\textwidth]{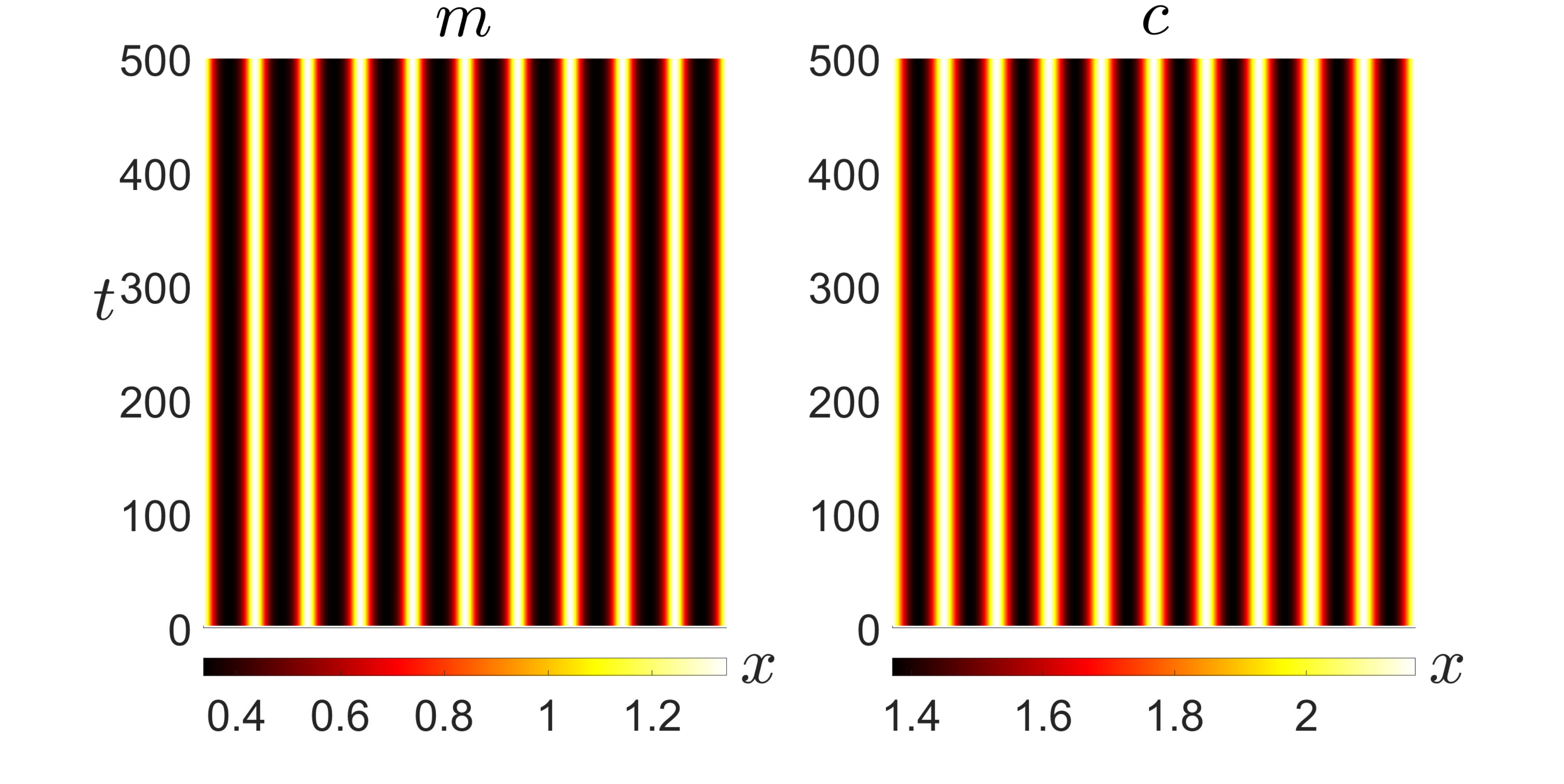}
\put(40,55){$\chi=3.5$}
\end{overpic}
}

\subfloat[M=0.02
\label{epsilon_0p08_M_0p02_chi_3p50_caso_2}]{
\begin{overpic}[width=0.35\textwidth]{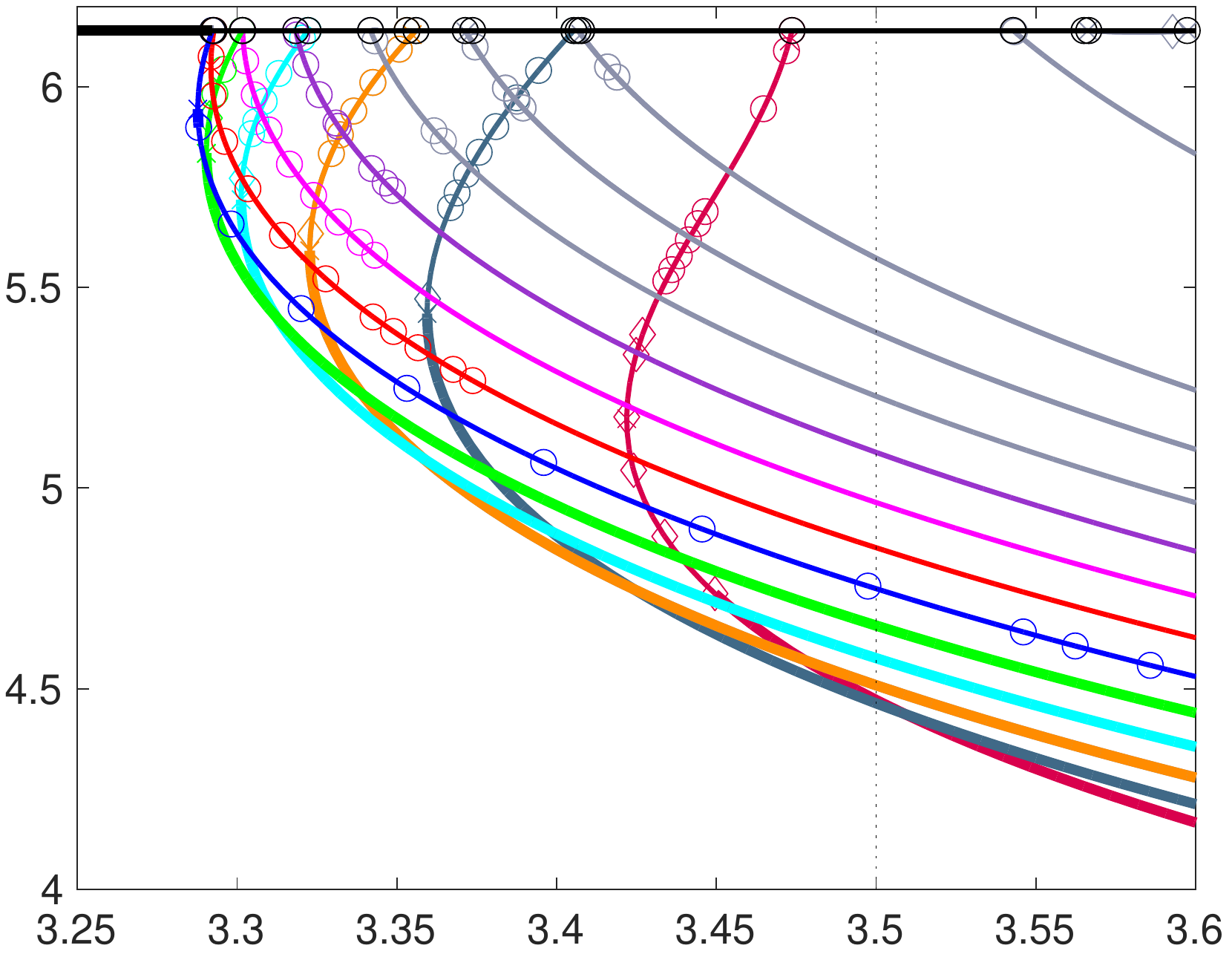}
\put(-20,50){$\|m\|_{L^2}$}
\put(50,-5){$\chi$}
\end{overpic}
\begin{overpic}[width=0.45\textwidth]{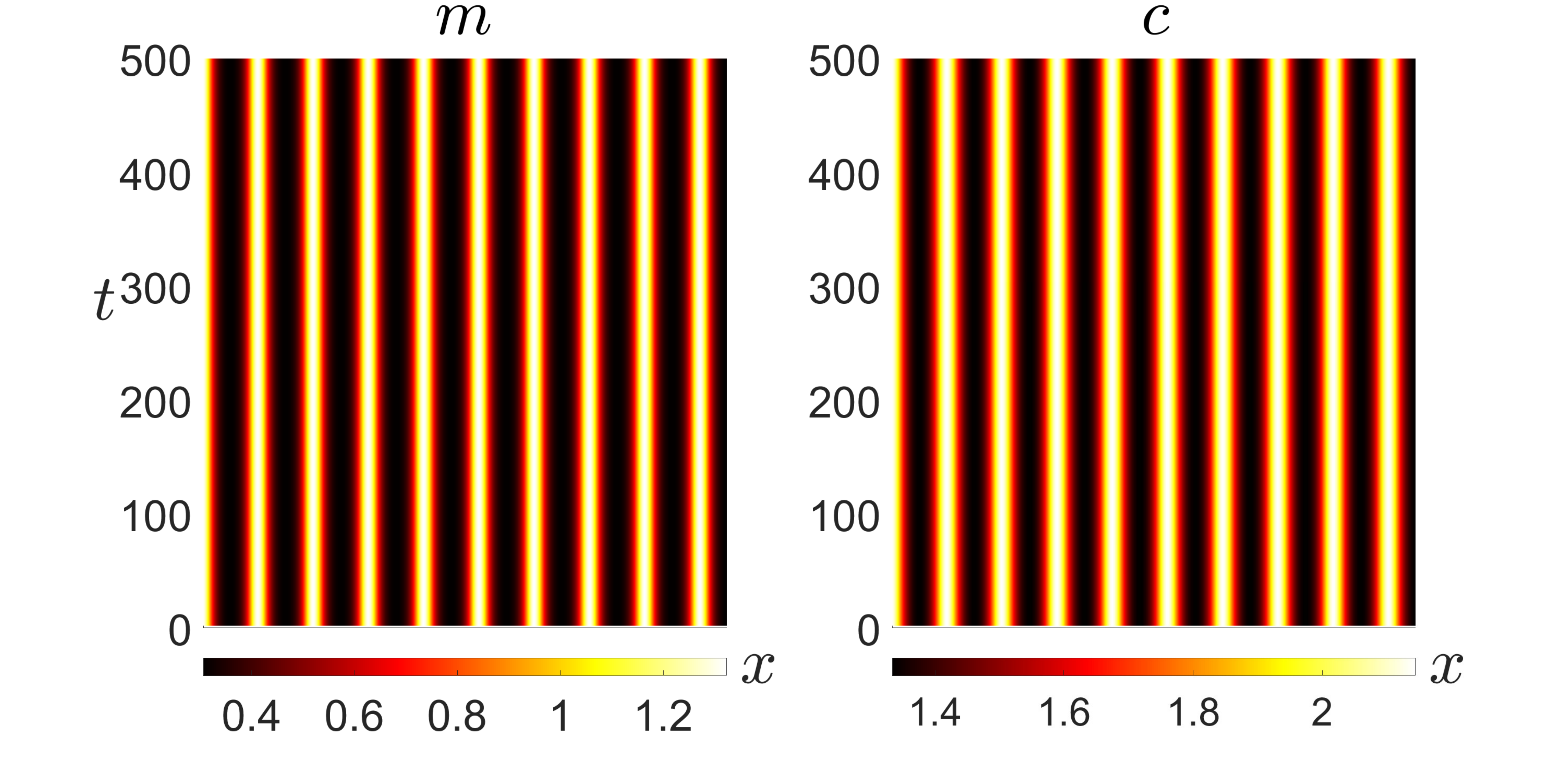}
\put(40,55){$\chi=3.5$}
\end{overpic}
}  
\caption{Influence of the parameter $M$ on the bifurcation diagram of steady-state solutions, obtained considering $\epsilon=0.08$ and $\Lambda$ as in Case 1. Dotted lines serve as references and denote the values for which numerical simulations (right panel) were performed. } 
\label{case1_M}
\end{figure}

In this Section, we investigate the influence of the Allee effect on stationary steady-state solutions, by combining numerical simulations and bifurcation diagrams. The numerical simulations are performed discretising system \eqref{nondim_syst} via finite differences in space and a first-order explicit method in time, starting from a small random perturbation of the homogeneous steady state. 
For the bifurcation diagram, we exploited the continuation software \texttt{pde2path}. Thanks to the results of the weakly nonlinear analysis, we can perform the numerical continuation on the reduced system with only two equations for $m$ and $c$, obtained from system \eqref{nondim_syst}, and considering the homogeneous solution $d=d_*$. In fact, the ODE part is stable and the remaining eigenvalues are solely determined by the PDEs part, i.e., the instability only arises from the PDEs subsystem. In the following bifurcation diagrams, thick lines denote stable solutions while thin lines are unstable ones. The homogeneous branch is shown in black, while branches bifurcating from the homogeneous one having the same number of peaks are indicated with the same colours when changing the parameters, to facilitate the comparison (for interpretation of the colour references, the reader is referred to the web version of this article). The continuation software \texttt{pde2path} runs on Matlab and the scripts needed for the numerical continuation of system \eqref{nondim_syst} are freely available in the GitHub folder \cite{GitHubFolder}. 
In the following, we want to explore the influence of the Allee effect on the outcomes of the model, thus we vary parameters $M$ and $\Lambda$. Parameter $\chi$ is taken as the bifurcation parameter, while we select two values for parameter $\epsilon$, namely $\epsilon=0.08$ and $\epsilon=0.8$. The former corresponds to a case in which the weakly nonlinear analysis up to the third order is able to predict the amplitude of the pattern close to the critical value $\chi_c$, while the latter is in the range of realistic values. All the other parameter values are taken as in \cite{lombardo2017demyelination} and reported here for convenience 
\begin{equation}\label{parameters}
\tau=1,\, \beta=1,\, r=1,\, \delta=1, \, \bar{L}=12\pi.  
\end{equation}
We first want to highlight the influence of the Allee effect with respect to the case of logistic growth investigated in \cite{lombardo2017demyelination}.
We fix $\epsilon=0.08$ and we vary the Allee parameter $M$. The parameter $\Lambda$ is taken as in Case 1, since this choice keeps the bifurcation points on the homogeneous branch fixed, allowing a direct comparison for selected values of the bifurcation parameter $\chi$. The results are shown in Figure \ref{case1_M}. In the left panel, the bifurcation diagrams with respect to the parameter $\chi$ are shown in the case of logistic growth (no Allee effect, for comparison) and with weak ($M=-0.5$ and $M=-0.1$) and strong ($M=0.02$) Allee effect. The homogeneous branch is shown in black, while branches bifurcating from the homogeneous corresponding to the same wavenumber (leading to a specific number of peaks in the solution profile) are indicated with the same colours when changing the parameters, to facilitate the comparison. Dotted vertical lines at $\chi=3.5$ are meant for reference and indicate the value used for the numerical simulations (right panel).
We notice that the bifurcation points on the homogeneous branch remain fixed, as expected. We observe major changes induced by an increase of parameter $M$:
\begin{itemize}
    \item[-] \textit{Loss of stability of the branches}. With an Allee effect and also increasing the parameter $M$, we observe a progressive loss of stability of the branches. For instance, the magenta branch is stable with a logistic growth but not with an Allee effect, and also the red and blue branches become unstable.
    \item[-] \textit{More prominent bend for subcritically bifurcating branches}. The branches bifurcate supercritically both with logistic growth and a weak Allee effect with $M=-0.5$ (qualitatively similar comparing the growth functions, see Figure \ref{fig:growth_cases1&2}). In particular, the first (blue) branch is stable for $\chi>\chi_{k^2}$, with $k^2=23$ (and corresponding solution presenting 11.5 peaks/bumps) while some of the others stabilise eventually. Increasing $M$, several branches become subcritical and present a wider stable region, namely for lower values of $\chi$. This is clearly visible for $M=0.02$. 
    \item[-] \textit{Stabilisation of spatial-periodic solutions with lower number of peaks}. For instance, setting $\chi=3.5$ and performing numerical simulations starting from a random perturbation of the homogeneous state, the numerical solution shows 11 peaks (blue branch) without Allee effect, 11.5 peaks (red branch) with $M=-0.5$, 10 peaks (cyan branch) with $M=-0.1$, and  9.5 peaks (orange branch) with $M=0.02$.
\end{itemize}

\begin{figure}[!ht]
\centering
\subfloat[$\Lambda$ in Case 1 \label{fig:epsilon_0p8_M_m0p5_chi_7p50_caso_2}]{
\begin{overpic}[width=0.4\textwidth]{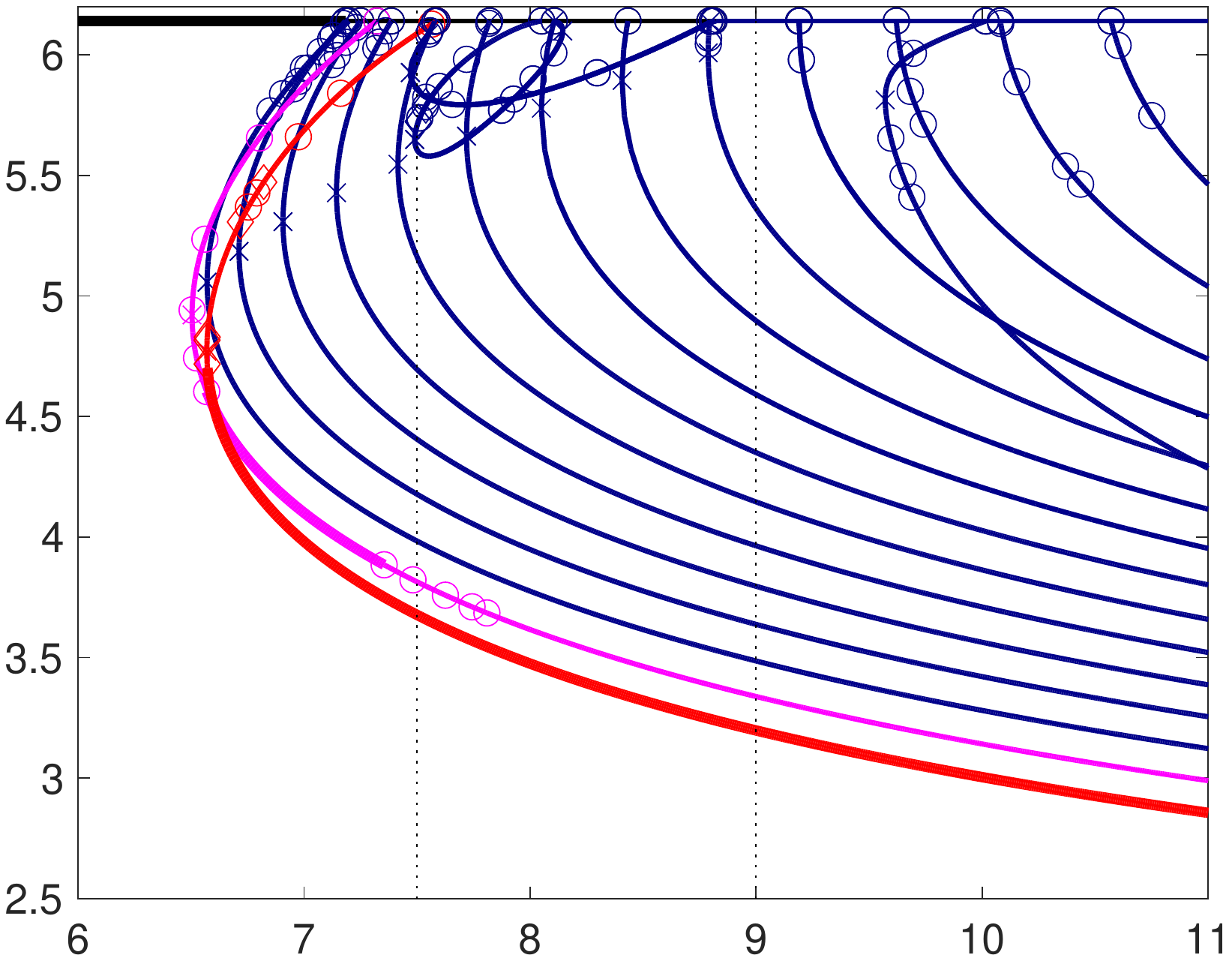}
\put(-20,50){$\|m\|_{L^2}$}
\put(50,-5){$\chi$}
\end{overpic}
\begin{overpic}[width=0.5\textwidth]{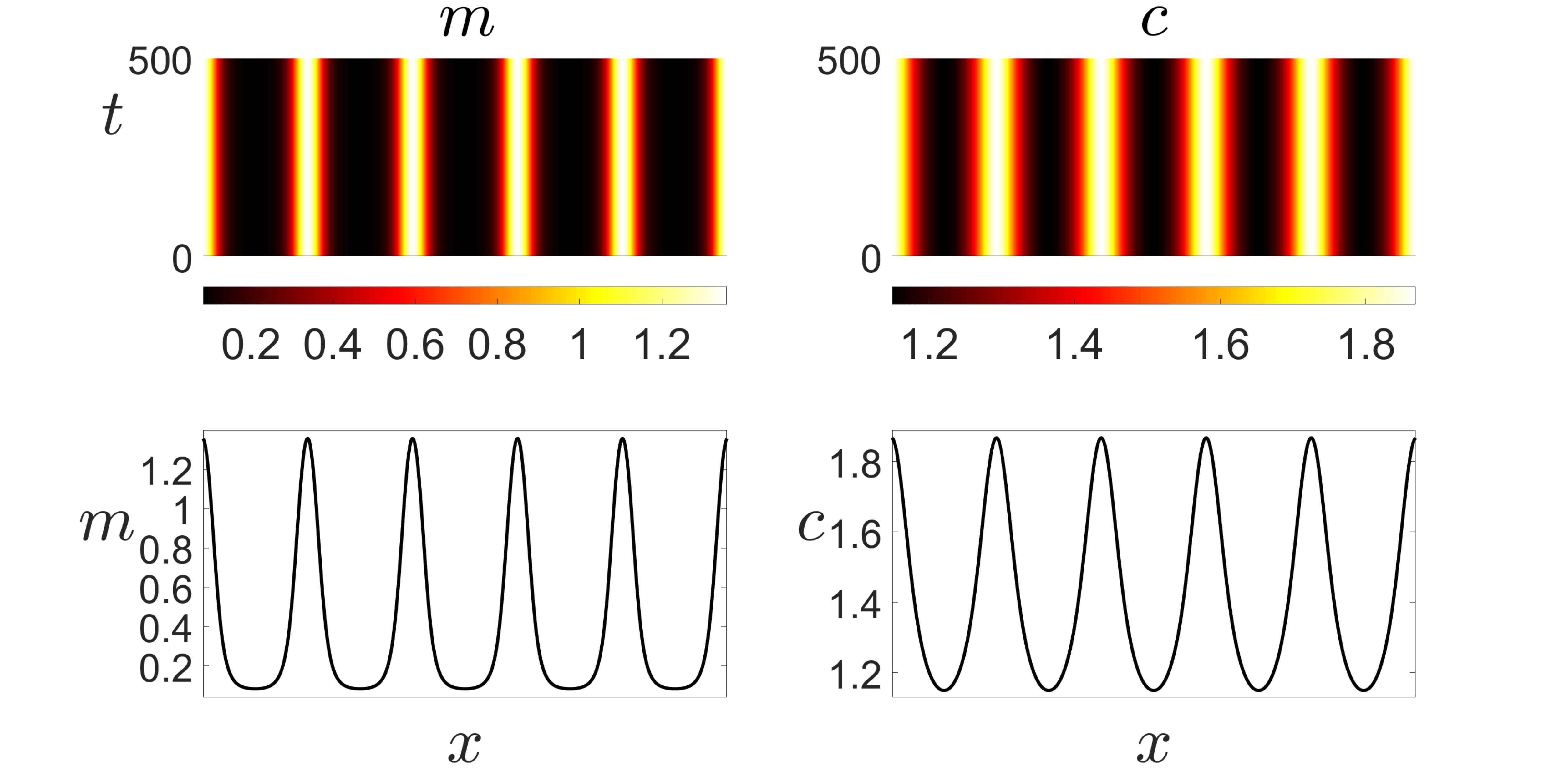}
\put(40,60){$\chi=7.5$}
\end{overpic}
}\\
\subfloat[$\Lambda$ in Case 2 \label{fig:epsilon_0p8_M_m0p5_chi_9_caso_1}]{
\begin{overpic}[width=0.4\textwidth]{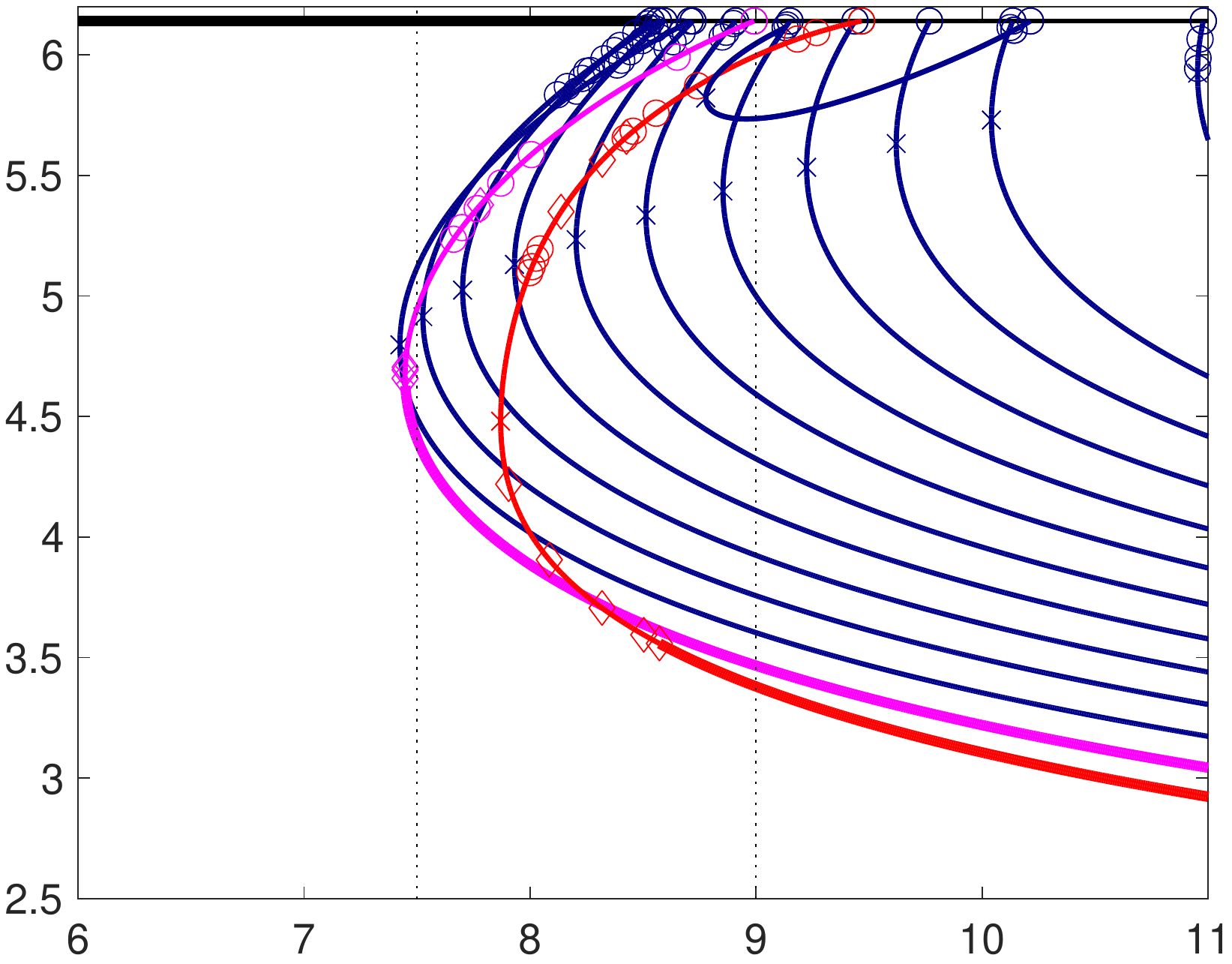}
\put(-20,50){$\|m\|_{L^2}$}
\put(50,-5){$\chi$}
\end{overpic}
\begin{overpic}[width=0.5\textwidth]{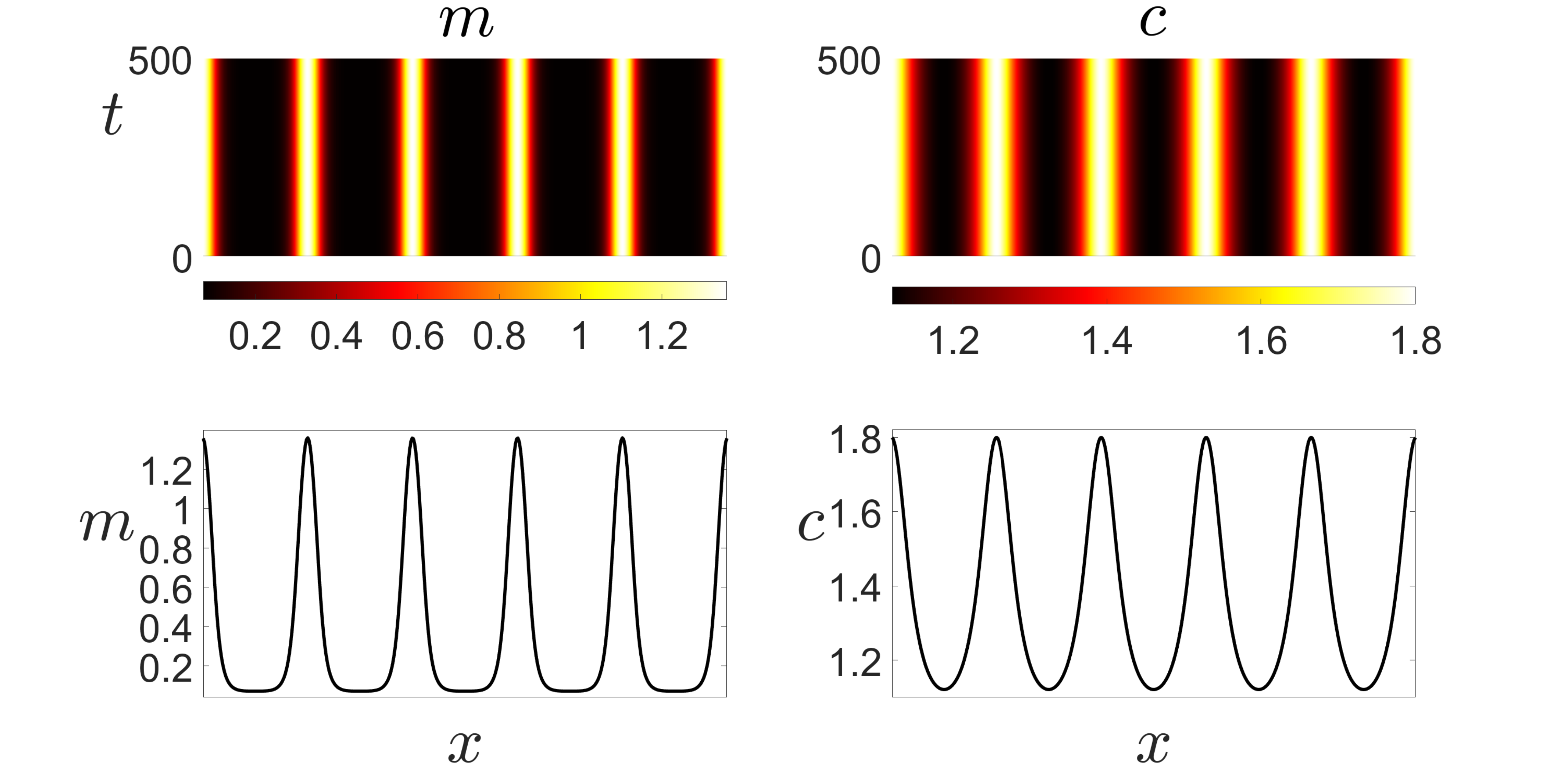}
\put(40,60){$\chi=9$}
\end{overpic}
}
\vspace{0.2cm}
\caption{Influence of the choice of parameter $\Lambda$ as in Case 1 (upper panel) and as in Case 2 (lower panel) on the bifurcation diagram of steady state solutions, obtained considering a realistic value $\epsilon=0.8$ and a weak Allee effect $M=-0.5$. In the bifurcation diagrams (left panel), the magenta branch corresponds to solutions with 5.5 peaks, while the red branch with 5 peaks. The homogeneous branch is denoted in black and all the other unstable branches are shown in blue. Dotted lines serve as references and denote the values for which numerical simulations (right panel) were performed.}\label{fig:case1VScase2}
\end{figure}

We now focus on the parameter $\Lambda$ and in particular on the difference between the two possible ways of defining it (Case 1 and Case 2 in Section \ref{sec:wna}). Remember that the choice of $\Lambda$ as in Case 2 allows maintaining the maximum growth rate as in the logistic growth (while in Case 1 we maintain the bifurcation point fixed, drastically reducing the maximum growth rate, see Figure \ref{fig:growth_cases1&2}). To investigate the influence of this choice on the bifurcation diagram and steady-state solutions, we consider $\epsilon=0.8$, corresponding to a realistic value for this parameter, and $M=-0.5$ (all the other parameter values are listed in \eqref{parameters}, except the bifurcation parameter $\chi$ which is specified in the text for the numerical simulations). Results are illustrated in Figure \ref{fig:case1VScase2}. In the bifurcation diagrams (left panel), the homogeneous branch is denoted in black and all the blue branches are unstable. Among the explored bifurcation points (with $\chi_{k^2}<12$), only two lead to regions of stable solutions and they are denoted in magenta and red. The magenta branch corresponds to solutions with 5.5 peaks, while the red branch with 5 peaks. The first difference comparing the bifurcation diagrams in Case 1 and Case 2 is, as expected, the position of the bifurcation points on the homogeneous branch: in particular, bifurcation points occur for lower values of $\chi$ in Case 1. The magenta and red branches are subcritical in both cases, but they are not qualitatively the same. At the two reference values $\chi=7.5$ and $\chi=9$, indicated by dotted lines, we have different stable solutions. In particular, at $\chi=7.5$, only the magenta branch turns out to be stable in Case 2, while in Case 1 only the red one is stable. On the contrary, at $\chi=9$, both the magenta and the red branch are stable in Case 2, while again in Case 1 only the red one is stable. Moreover, we notice that to have a similar steady state solution (right panel), we have to select two different values of $\chi$. In fact, we show in the right panel the numerical simulations obtained for $\chi=9$ in Case 2 and $\chi=7.5$ in Case 1 and the profile of the steady state solutions. Thus, the choice of $\Lambda$ seems to have a qualitative and quantitative impact on the possible outcome of the system, which is crucial in applications. 

\clearpage
\section{Conclusion and outlook}\label{sec:conclusion}

In this paper, we proposed a modification of the mathematical model describing inflammation and demyelination patterns in the brain caused by Multiple Sclerosis proposed in \cite{lombardo2017demyelination}. In particular, we hypothesized that a minimal amount of macrophages is needed to be able to start and sustain the inflammatory response, while if the concentration of macrophages is too low, then there is no inflammatory response. To model this, we introduced an Allee effect in the model function for macrophage activation. We investigated the emergence and the type of Turing patterns on a 1D domain by combining linearised and weakly nonlinear analysis, bifurcation diagrams and numerical simulations, focusing on the comparison with the previous model \cite{lombardo2017demyelination}, formulated using a logistic growth for macrophage activation. We selected two cases for this comparison, one in which the bifurcation points on the homogeneous branch are the same in both the logistic and Allee case, the other in which the maximum macrophage activation rates are the same (even if it is attained for different values, see Figure \ref{fig:growth_cases1&2}).

When an Allee effect is considered, we observed clear trends as the Allee parameter $M$ increases. In detail, branches which are stable in the logistic case become unstable, leading to very few stable branches in the bifurcation diagrams. We also observed that the Allee effect leads more often to subcritical bifurcations or to a more prominent bend. Finally, we also observed the stabilisation of branches bifurcating for greater values of the bifurcation parameter (very far from the first bifurcation point and the critical value), leading to stable solutions with a lower number of peaks. 

We plan to further investigate this model from several points of view. 

In cases when the Landau coefficient $L$ is negative and thus the third order Stuart--Landau equation is not able to capture the amplitude of the pattern, one should push the weakly nonlinear analysis to a higher order, obtaining a quintic Stuart--Landau equation for the amplitude of the pattern as in \cite{lombardo2017demyelination}. 

Another phenomenon worth to be investigated is the wavefront invasion, namely the existence of travelling wavefronts connecting two different steady solutions of the equations; this would require a weakly nonlinear analysis similar to that performed in Section \ref{sec:wna}, but with an additional expansion of the derivative with respect to the space variable, in order to take into account the slow modulation in space of pattern amplitude \cite{barresi2016wavefront}.

The analysis and the numerical simulations in this paper are performed on a one-dimensional domain. We plan to extend the results to two-dimensional domains, in order to check the formation of different types of spatial patterns in presence of Allee effect varying the parameters and to compare the results with the ones in \cite{lombardo2017demyelination} for the model without Allee effect. On the other hand, the presence of time-periodic spatial patterns was already noted in \cite{lombardo2017demyelination}. The effect of the Allee parameters on these patterns will also be a matter of future work.


Moreover, since we observed that the Allee effect sharpens the subcritical bifurcations, the bifurcation structure might show snaking branches leading to localised patterns. In particular, snaking branches might arise from secondary bifurcation points on branches bifurcating \textit{highly} subcritically at $\chi_{k^2}\gg\chi_c$. To identify such branches,  the same weakly nonlinear analysis can be extended at the general bifurcation value $\chi_{bif}$ (with similar but slightly different computations) and looking at the sign and value of the Landau coefficient~$L$.\\[0.3cm]

\textbf{Acknowledgments.}
The authors are grateful to the Editors for the invitation to contribute to this special issue to honour the memory of Prof.~Salvatore Rionero, an outstanding scientist whose findings on stability and bifurcations in reaction-diffusion systems also inspired their steps in this direction.\\

\textbf{Funding.} 
The authors are members of the INdAM-GNFM Group. The work has been performed in the frame of the Italian National Research Projects ``Multiscale phenomena in Continuum Mechanics: singular limits, off-equilibrium and transitions'' (Prin 2017YBKNCE) and ``Integrated Mathematical Approaches to Socio--Epidemiological Dynamics'' (Prin 2020JLWP23, CUP: E15F21005420006).
GM thanks the support of INdAM-GNFM with the research project ``From kinetic to macroscopic models for tumor--immune system competition''.
CS has been partially supported by the Fondazione ACRI in the framework of the Research Awards Young Investigator Training Programme 2019. MB, MG, CS also thank the support of the University of Parma with the FIL project ``Collective and self-organised dynamics: kinetic and network approaches''.
\bibliographystyle{plain}
\bibliography{bibliography}


\end{document}